\documentclass[10pt]{article}

\usepackage{amscd,amsmath, amssymb, fancyhdr, epsfig}
\usepackage[matrix,arrow]{xy}

\numberwithin{equation}{section}

\newcommand{\version}{version 4.0,\ \   Dec. 01, 2008}

\def\eqref#1{(\ref{#1})}
\newcommand{\goth}{\mathfrak}

\newcommand{\arrow}{{\:\longrightarrow\:}}

\newcommand{\C}{{\Bbb C}}
\newcommand{\R}{{\Bbb R}}

\newcommand{\6}{\partial}
\def\1{\sqrt{-1}\:}
\newcommand{\restrict}[1]{{\left|_{{\phantom{|}\!\!}_{#1}}\right.}}
\newcommand{\cntrct}                
{\hspace{2pt}\raisebox{1pt}{\text{$\lrcorner$}}\hspace{2pt}}

\makeatletter
\def\x@arrow{\DOTSB\Relbar}
\def\xlongequalsignfill@{\arrowfill@\x@arrow\Relbar\x@arrow}
\newcommand{\xlongequal}[2][]{%
        \ext@arrow 0099\xlongequalsignfill@{#1}{#2}}
\def\xlongrightarrowfill@{\arrowfill@\relbar\relbar\longrightarrow}
\newcommand{\xlongrightarrow}[2][]{%
	\ext@arrow 0099\xlongrightarrowfill@{#1}{#2}}
\makeatother

\newcommand{\calo}{{\cal O}}

\renewcommand{\bar}{\overline}
\renewcommand{\phi}{\varphi}
\renewcommand{\epsilon}{\varepsilon}
\renewcommand{\geq}{\geqslant}
\renewcommand{\leq}{\leqslant}

\newcommand{\End}{\operatorname{End}}

\newcommand{\Id}{\operatorname{Id}}
\newcommand{\Cl}{\operatorname{Cl}}

\newcommand{\const}{\operatorname{\text{\sf const}}}
\newcommand{\Vol}{\operatorname{Vol}}

\newcommand{\Hol}{\operatorname{Hol}}

\newcommand{\Alt}{\operatorname{Alt}}

\newcommand{\Spin}{\operatorname{Spin}}

\newcounter{Mycounter}[section]
\newcounter{lemma}[section]
\renewcommand{\thelemma}{{Lemma \thesection.\arabic{lemma}}}
\newcommand{\lemma}{%
    \setcounter{lemma}{\value{Mycounter}}
    \refstepcounter{lemma}
    \stepcounter{Mycounter}
    {\noindent \bf \thelemma:\ }}

\newcounter{claim}[section]
\newcounter{sublemma}[section]
\newcounter{corollary}[section]
\renewcommand{\thecorollary}{{Corollary \thesection.\arabic{corollary}}}
\newcommand{\corollary}{%
    \setcounter{corollary}{\value{Mycounter}}
    \refstepcounter{corollary}
    \stepcounter{Mycounter}
    {\noindent \bf \thecorollary:\ }}

\newcounter{theorem}[section]
\renewcommand{\thetheorem}{{Theorem \thesection.\arabic{theorem}}}
\newcommand{\theorem}{%
    \setcounter{theorem}{\value{Mycounter}}
    \refstepcounter{theorem}
    \stepcounter{Mycounter}
    {\noindent \bf \thetheorem:\ }}

\newcounter{conjecture}[section]
\renewcommand{\theconjecture}{{Conjecture \thesection.\arabic{conjecture}}}
\newcommand{\conjecture}{%
    \setcounter{conjecture}{\value{Mycounter}}
    \refstepcounter{conjecture}
    \stepcounter{Mycounter}
    {\noindent \bf \theconjecture:\ }}

\newcounter{proposition}[section]
\renewcommand{\theproposition}
      {{Proposition \thesection.\arabic{proposition}}}
\newcommand{\proposition}{%
    \setcounter{proposition}{\value{Mycounter}}
    \refstepcounter{proposition}
    \stepcounter{Mycounter}
    {\noindent \bf \theproposition:\ }}

\newcounter{definition}[section]
\renewcommand{\thedefinition}
      {{Definition~\thesection.\arabic{definition}}}
\newcommand{\definition}{%
    \setcounter{definition}{\value{Mycounter}}
    \refstepcounter{definition}
    \stepcounter{Mycounter}
    {\noindent \bf \thedefinition:\ }}

\newcounter{example}[section]
\newcounter{remark}[section]
\renewcommand{\theremark}{{Remark \thesection.\arabic{remark}}}
\newcommand{\remark}{%
    \setcounter{remark}{\value{Mycounter}}
    \refstepcounter{remark}
    \stepcounter{Mycounter}
    {\noindent \bf \theremark:\ }}

\newcounter{problem}[section]
\newcounter{question}[section]
\makeatletter

\@addtoreset{equation}{section} \@addtoreset{footnote}{section}
\makeatother

\def\blacksquare{\hbox{\vrule width 5pt height 5pt depth 0pt}}
\def\endproof{\blacksquare}

\begin{document}
\begin{center}
{\LARGE\bf
Balanced HKT metrics and strong HKT\\[1mm] metrics on
hypercomplex manifolds\\[3mm]
}

 Misha
Verbitsky\footnote{Misha Verbitsky is 
supported by CRDF grant RM1-2354-MO02.}

\end{center}

{\small \hspace{0.10\linewidth}
\begin{minipage}[t]{0.85\linewidth}
{\bf Abstract} \\ A manifold $(M,I,J,K)$ is called
hypercomplex if $I,J,K$ are complex structures satisfying
quaternionic relations. A quaternionic Hermitian 
hypercomplex manifold
is called HKT (hyperk\"ahler with torsion) if the
(2,0)-form $\Omega$ associated with the corresponding
$Sp(n)$-structure satisfies $\6\Omega=0$.
A Hermitian metric $\omega$ on a complex manifold is 
called balanced if $d^*\omega=0$. We show that
balanced HKT metrics are precisely the quaternionic 
Calabi-Yau metrics defined 
in terms of the quaternionic Monge-Amp\`ere equation.
In particular, a balanced HKT-metric is unique
in its cohomology class, and it always exists
if the quaternionic Calabi-Yau theorem is true.
We investigate the cohomological properties of 
strong HKT metrics (the quaternionic Hermitian 
metrics, satisfying, in addition to the HKT condition,
the relation $dd^c \omega=0$), and show that
the space of strong HKT metrics is finite-dimensional.
Using Howe's duality for representations
of $Sp(n)$, we prove a hyperk\"ahler version of 
Hodge-Riemann bilinear relations. We use it to show that 
a manifold admitting a balanced HKT-metric
never admits a strong HKT-metric, if $\dim_\R M \geq 12$.
\end{minipage}
}

\tableofcontents


\section{Introduction}


The notion of an HKT manifold 
was introduced by the physicists,
but it proved to be immensely useful in matematics.

A {\bf hypercomplex manifold} is a manifold equipped with
almost complex structure operators $I, J, K:\; TM \arrow TM$,
integrable and satisfying the standard quaternionic relations
$I^2=J^2=K^2= IJK = -\Id_{TM}$.

This gives a quaternionic algebra action on $TM$;
the group $Sp(1) \cong SU(2)$ of unitary quaternions
acts on all tensor powers of $TM$ by multiplicativity.

A {\bf quaternionic Hermitian structure}
 on a hypercomplex manifold is an $SU(2)$-invariant Riemannian
metric. Such a metric gives a reduction of the structure
group of $M$ to $Sp(n)= U(n, {\Bbb H})$. 

With any quaternionic Hermitian structure on $M$
one associates a non-degenerate $(2,0)$-form
$\Omega \in \Lambda^{2,0}_I(M)$, as 
follows.\footnote{$\Lambda^*(M)$ denotes the bundle
of differential forms, and $\Lambda^*(M)=\oplus_{p,q}\Lambda^{p,q}_I(M)$ 
its Hodge decomposition, taken with respect to the
complex structure $I$ on $M$.}
Consider the differential forms
\begin{equation}\label{_three_2_forms_qH_Equation_}
\omega_I(\cdot, \cdot) := g(\cdot, I\cdot), \ \ \omega_J(\cdot, \cdot) := g(\cdot, J\cdot), \ \ \omega_K(\cdot, \cdot) := g(\cdot, K\cdot).
\end{equation}
It is easy to check that the form $\Omega:= \omega_J + \1\Omega_K$
is of Hodge type $(2,0)$ with respect to $I$. 

If the form $\Omega$ is closed, one has
$d\omega_I = d \omega_J = d\omega_K =0$, and the
manifold $(M,I,J,K, g)$ is called 
{\bf hyperk\"ahler} (\cite{_Besse:Einst_Manifo_}). 
The hyperk\"ahler condition is very restrictive.

A hypercomplex, quaternionic Hermitian manifold
$(M,I,J,K, g)$ is called {\bf an HKT-manifold} 
(hyperk\"ahler with torsion) if $\6\Omega=0$, where 
$\6$ denotes the $(1,0)$-part of the differential. 
In other words, a manifold is HKT
if $d\Omega\in \Lambda_I^{2,1}(M)$. 

The form $\Omega \in \Lambda_I^{2,0}(M)$
is called {\bf an HKT-form on $(M,I,J,K)$}. 

\hfill

\remark\label{_Herm_via_(2,0)_Remark_}
The quaternionic Hermitian form $g$ can be easily
reconstructed from $\Omega$. Indeed, for any
$x, y \in T^{1,0}_I(M)$, one has
\[
2 g(x, \bar y)= \Omega(x, J(\bar y)),
\]
as a trivial calculation implies.

\hfill

HKT-manifolds were first introduced by the physicists
(\cite{_Howe_Papado_}; see also \cite{_Gra_Poon_})
in a completely different context. Given a complex Hermitian
manifold $(M, I, g)$, one defines {\bf a Bismut connection}
$\nabla:\; TM \arrow TM \otimes \Lambda^1M$, 
determined by the following properties
\begin{description}
\item[(i)] $\nabla I = \nabla g =0$
\item[(ii)] The torsion form $T_\nabla \in TM \otimes \Lambda^2M$
is totally antisymmetric, if one identifies $TM \otimes \Lambda^2M$
with $\Lambda^1 M \otimes \Lambda^2M$.
\end{description}
This connection has its origins also in physics,
due to A. Strominger (\cite{_Strominger:Bismut_}),
who defined it earlier than Bismut's paper appeared.
It is well known such $\nabla$ exists, and it is unique
(see e.g. \cite{_Gauduchon:Herm+Dirac_}). 
The torsion 3-form $T$ of Bismut connection
can be written down explicitly in terms
of its Hermitian form $\omega$: \[ T= - I d\omega.\]

Now, suppose that $(M,I,J,K, g)$ is a hypercomplex,
quaternionic Hermitian manifold. The metric $g$ can
be used to define the Bismut connections 
$\nabla_I, \nabla_J, \nabla_K$ associated with $I, J, K$.
It is known (see e.g. \cite{_Gra_Poon_}) that
$\nabla_I = \nabla_J = \nabla_K$ if and only if
$(M,I,J,K, g)$ is HKT. This was the original definition
of HKT structures (\cite{_Howe_Papado_}).

\hfill

\remark \label{_HKT_via_torsion_Remark_}
Let $(M, g)$ be a Riemannian manifold,
and $\nabla$ a connection on $M$ which satisfies $\nabla g=0$.
Such a connection is uniquely determined by its torsion
form; this is proven by the same argument as used
to show existence and uniqueness of the Levi-Civita
connection. However, the torsion term of
the Bismut connection is written as $T= - I d\omega$.
Therefore,  $\nabla_I = \nabla_J = \nabla_K$
is equivalent to the following relation:
\[
 - I d\omega_I =  - J d\omega_J =  - K d\omega_K
\]
This relation can be used as one more definition
of HKT structures.


\section{Quaternionic Dolbeault complex on a hypercomplex manifold}


\subsection{Quaternionic Dolbeault complex: a definition}
\label{_qD_Subsection_}

It is well-known that any irreducible representation
of $SU(2)$ over $\C$ can be obtained as a symmetric power
$S^i(V_1)$, where $V_1$ is a fundamental 2-dimensional
representation. We say that a representation $W$ 
{\bf has weight $i$} if it is isomorphic to $S^i(V_1)$.
A representation is said to be {\bf pure of weight $i$}
if all its irreducible components have weight $i$.

\remark\label{_weight_multi_Remark_}
The Clebsch-Gordan formula (see \cite{_Humphreys_})
claims that the weight is {\em multiplicative}, 
in the following sense: if $i\leq j$, then
\[
V_i\otimes V_j = \bigoplus_{k=0}^i V_{i+j-2k},
\]
where $V_i=S^i(V_1)$ denotes the irreducible
representation of weight $i$.

\hfill

Let $M$ be a hypercomplex  manifold,
$\dim_{\Bbb H}M=n$.
There is a natural multiplicative action of $SU(2)\subset
{\Bbb H}^*$ on $\Lambda^*(M)$, associated with the
hypercomplex structure. 

\hfill

Let $V^i\subset \Lambda^i(M)$ be a maximal
$SU(2)$-invariant subspace of weight $<i$.
The space $V^i$ is well defined, because
it is a sum of all irreducible representations
$W\subset \Lambda^i(M)$ of weight $<i$.
Since the weight is multiplicative
(\ref{_weight_multi_Remark_}), $V^*= \bigoplus_i V^i$
is an ideal in $\Lambda^*(M)$. 

It is easy to see that the de Rham differential
$d$ increases the weight by 1 at most. Therefore,
$dV^i\subset V^{i+1}$, and $V^*\subset \Lambda^*(M)$
is a differential ideal in the de Rham DG-algebra
$(\Lambda^*(M), d)$.

\hfill

\definition\label{_qD_Definition_}
Denote by $(\Lambda^*_+(M), d_+)$ the quotient algebra
$\Lambda^*(M)/V^*$
It is called {\bf the quaternionic Dolbeault algebra of
  $M$}, or {\bf the quaternionic Dolbeault complex} 
(qD-algebra or qD-complex for short).

\hfill

\remark 
The complex $(\Lambda^*_+(M), d_+)$ 
was constructed much earlier by Capria and Salamon,
(\cite{_Capria-Salamon_}) in a different (and much 
more general) situation, and much studied since then.

\subsection{The Hodge decomposition of the quaternionic
  Dolbeault complex}.
\label{_Hodge_on_qD_Subsection_}

The Hodge bigrading is compatible with the weight decomposition
of $\Lambda^*(M)$, and gives a Hodge decomposition 
of $\Lambda^*_+(M)$ (\cite{_Verbitsky:HKT_}):
\[
\Lambda^i_+(M) = \bigoplus_{p+q=i}\Lambda^{p,q}_{+,I}(M).
\]
The spaces $\Lambda^{p,q}_{+,I}(M)$
are the weight spaces for a particular choice of a Cartan
subalgebra in $\goth{su}(2)$. The $\goth{su}(2)$-action
induces an isomorphism of the weight spaces
within an irreducible representation. This
gives the following result.

\hfill

\proposition \label{_qD_decompo_expli_Proposition_}
Let $(M,I,J,K)$ be a hypercomplex manifold and
\[
\Lambda^i_+(M) = \bigoplus_{p+q=i}\Lambda^{p,q}_{+,I}(M)
\]
the Hodge decomposition of qD-complex defined above.
Then there is a natural isomorphism
\begin{equation}\label{_qD_decompo_Equation_}
\Lambda^{p,q}_{+,I}(M)\cong \Lambda^{p+q,0}(M,I).
\end{equation}

{\bf Proof:} See \cite{_Verbitsky:HKT_}. \endproof

\hfill

This isomorphism is compatible with a natural algebraic
structure on \[ \bigoplus_{p+q=i}\Lambda^{p+q,0}(M,I),\]
and with the Dolbeault differentials, in the following
way.

\hfill

Let $(M,I,J,K)$ be a hypercomplex manifold.
We extend \[ J:\; \Lambda^1(M) \arrow \Lambda^1(M)\]
to $\Lambda^*(M)$ by multiplicativity. Recall that 
\[ J(\Lambda^{p,q}(M,I))=\Lambda^{q,p}(M,I), \]
because $I$ and $J$ anticommute on $\Lambda^1(M)$.
Denote by 
\[ \6_J:\;  \Lambda^{p,q}(M,I)\arrow \Lambda^{p+1,q}(M,I)
\]
the operator $J\circ \bar\6 \circ J$, where
$\bar\6:\;  \Lambda^{p,q}(M,I)\arrow \Lambda^{p,q+1}(M,I)$
is the standard Dolbeault operator on $(M,I)$, that is, the
$(0.1)$-part of the de Rham differential.
Since $\bar\6^2=0$, we have $\6_J^2=0$.
In \cite{_Verbitsky:HKT_} it was shown that $\6$ and $\6_J$
anticommute:
\begin{equation}\label{_commute_6_J_6_Equation_}
\{\6_J, \6 \}=0.
\end{equation}

Consider the quaternionic Dolbeault complex
$(\Lambda^*_+(M), d_+)$ constructed in Subsection
\ref{_qD_Subsection_}. Using the Hodge bigrading, we can
decompose this complex, obtaining a bicomplex 
\[
\Lambda^{*, *}_{+,I}(M) \xlongrightarrow{d^{1,0}_{+,I}, d^{0,1}_{+,I}}
\Lambda^{*, *}_{+,I}(M)
\]
where $d^{1,0}_{+,I}$,  $d^{0,1}_{+,I}$ are the Hodge components of 
the quaternionic Dolbeault differential $d_+$, taken with
respect to $I$. 

\hfill

\theorem\label{_bico_ide_Theorem_}
Under the multiplicative isomorphism 
\[
\Lambda^{p,q}_{+,I}(M)\cong \Lambda^{p+q,0}(M,I)
\]
constructed in \ref{_qD_decompo_expli_Proposition_},
$d^{1,0}_+$
corresponds to $\6$ and $d^{0,1}_+$
to $\6_J$:

\begin{equation}\label{_bicomple_XY_Equation}
\begin{minipage}[m]{0.85\linewidth}
{\tiny $
\xymatrix @C+1mm @R+10mm@!0  { 
  && \Lambda^0_+(M) \ar[dl]^{d^{0,1}_+} \ar[dr]^{d^{1,0}_+} 
   && && && \Lambda^{0,0}_I(M) \ar[dl]^{\6} \ar[dr]^{ \6_J}
   &&  \\
 & \Lambda^{1,0}_+(M) \ar[dl]^{d^{0,1}_+} \ar[dr]^{d^{1,0}_+} &
 & \Lambda^{0,1}_+(M) \ar[dl]^{d^{0,1}_+} \ar[dr]^{d^{1,0}_+}&& 
\text{\large $\cong$} &
 &\Lambda^{1,0}_I(M)\ar[dl]^{ \6} \ar[dr]^{ \6_J}&  &
 \Lambda^{1,0}_I(M)\ar[dl]^{ \6} \ar[dr]^{ \6_J}&\\
 \Lambda^{2,0}_+(M) && \Lambda^{1,1}_+(M) 
   && \Lambda^{0,2}_+(M)& \ \ \ \ \ \ & \Lambda^{2,0}_I(M)& & 
\Lambda^{2,0}_I(M) & &\Lambda^{2,0}_I(M) \\
}
$
}
\end{minipage}
\end{equation}
Moreover, under this isomorphism, 
$\omega_I\in \Lambda^{1,1}_{+,I}(M)$ corresponds to
$\Omega\in\Lambda^{2,0}_I(M)$. 

{\bf Proof:} See \cite{_Verbitsky:HKT_} or \cite{_Verbitsky:qD_}.
 \endproof

\subsection{Positive $(2,0)$-forms on hypercomplex
  manifolds}
\label{_posi_2,0-forms_Subsection_}

The notion of positive $(2p,0)$-forms on hypercomplex
manifolds (sometimes called q-positive, or ${\Bbb H}$-positive)
was developed in \cite{_V:reflexive_} and 
\cite{_Alesker_Verbitsky_HKT_} (see also
\cite{_AV:Calabi_}
and \cite{_Verbitsky:skoda.tex_}).
For our present purposes, only $(2,0)$-forms are interesting,
but everything can be immediately generalized to a
general situation

\hfill

Let $\eta\in \Lambda^{p,q}_I(M)$ be a 
differential form. Since $I$ and $J$ anticommute,
$J(\eta)$ lies in $\Lambda^{q,p}_I(M)$.
Clearly, $J^2\restrict {\Lambda^{p,q}_I(M)}=(-1)^{p+q}$.
For $p+q$ even, $J\restrict {\Lambda^{p,q}_I(M)}$
is an anticomplex involution, that is, a real structure
on $\Lambda^{p,q}_I(M)$. 
A form 
$\eta \in \Lambda^{2p,0}_I(M)$ is called {\bf real} if
$J(\bar\eta)=\eta$. 

For a real $(2,0)$-form $\eta$, 
\[ 
   \eta\left(x, J(\bar x))\right)=
   \bar \eta\left(J(x), J^2 (\bar x)\right)=
 \bar \eta\left(\bar x, J(x)\right),
\] 
for any $x \in T^{1,0}_I(M)$.
From a definition of a real form,
we obtain that the scalar $\eta\left(x, J(\bar x)\right)$
is always real.

\hfill

\definition
A real $(2,0)$-form $\eta$ on a hypercomplex manifold
is called {\bf positive} if 
$\eta\left(x, J(\bar x)\right)\geq 0$ 
for any $x \in T^{1,0}_I(M)$, and {\bf strictly positive}
if this inequality is strict, for all $x\neq 0$.

\hfill

An HKT-form $\Omega\in \Lambda^{2,0}_I(M)$ of any
HKT-structure is strictly positive, as follows from 
\ref{_Herm_via_(2,0)_Remark_}. Moreover, HKT-structures
on a hypercomplex manifold are in one-to-one
correspondence with closed, strictly positive 
$(2,0)$-forms.

The analogy between K\"ahler forms and HKT-forms can be
pushed further; it turns out that any 
HKT-form $\Omega\in \Lambda^{2,0}_I(M)$ has a local
potential $\phi\in C^\infty(M)$, in such a way
that $\6\6_J\phi=\Omega$ (\cite{_Alesker_Verbitsky_HKT_}). 
Here $\6\6_J$ is a composition of $\6$ and $\6_J$
defined on quaternionic Dolbeault complex as above
(these operators anticommute).


\section{$SL(n, {\Bbb H})$-manifolds}


\subsection{An introduction to $SL(n, {\Bbb H})$-geometry}

As Obata has shown (\cite{_Obata_}), a hypercomplex manifold
$(M,I,J,K)$ admits a necessarily unique torsion-free connection, preserving
$I,J,K$. The converse is also true: if a manifold 
$M$ equipped with an action of ${\Bbb H}$
admits a torsion-free connection preserving
the quaternionic action, it is hypercomplex. This implies that a
hypercomplex structure on a manifold can be defined as a
torsion-free connection with holonomy in $GL(n, {\Bbb  H})$.
This connection is called {\bf the Obata connection}
on a hypercomplex manifold. 

Connections with restricted holonomy are one of the central notions
in Riemannian geometry, due to Berger's classification of
irreducible holonomy of Riemannian manifolds. However, a similar
classification exists for a general torsion-free connection
(\cite{_Merkulov_Sch:long_}). In the Merkulov-Schwachh\"ofer list,
only three subroups of $GL(n, {\Bbb H})$ occur. In addition to the
compact group $Sp(n)$ (which defines hyperk\"ahler geometry), also
$GL(n, {\Bbb H})$ and its commutator $SL(n, {\Bbb H})$ appear,
corresponding to hypercomplex manifolds and hypercomplex manifolds
with trivial determinant bundle, respectively. Both of these
geometries are interesting, rich in structure and examples, and
deserve detailed study.

It is easy to
see that $(M,I)$ has holomorphically trivial canonical
bundle, for any $SL(n, {\Bbb H})$-manifold 
$(M,I, J, K)$ (\cite{_Verbitsky:canoni_}). For a
hypercomplex manifold with trivial canonical bundle
admitting an HKT metric, a version of Hodge
theory was constructed (\cite{_Verbitsky:HKT_}). Using this 
result, it was shown that a compact hypercomplex manifold with trivial
canonical bundle has holonomy in $SL(n,{\Bbb H})$, if it admits an
HKT-structure (\cite{_Verbitsky:canoni_}).

In \cite{_BDV:nilmanifolds_}, it was shown that holonomy of all
hypercomplex nilmanifolds lies in $SL(n, {\Bbb H})$.
Many (probably, most) working examples of hypercomplex manifolds are
in fact nilmanifolds, and by this result they all belong to the
class of $SL(n, {\Bbb H})$-manifolds. 

The $SL(n, {\Bbb H})$-manifolds were studied in \cite{_AV:Calabi_}  
and \cite{_Verbitsky:skoda.tex_}, 
because on such manifolds the quaternionic Dolbeault complex 
is identified with a part of de Rham complex (\ref{_V_main_Proposition_}).
 Under this identification, ${\Bbb H}$-positive forms become positive
in the usual sense, and $\6$, $\6_J$-closed or exact forms
become $\6, \bar\6$-closed or exact. This linear-algebraic
identification is especially useful in the study of
quaternionic Monge-Amp\`ere equation (Subsection 
\ref{_qCY_Subsection_}).

\subsection{The map ${\cal V}_{p,q}:\;
  \Lambda^{p+q,0}_I(M)\arrow\Lambda^{n+p, n+q}_I(M)$\\
on $SL(n, {\Bbb H})$-manifolds}
\label{_V_p,q_Subsection_}

Let $(M,I,J,K)$ be an $SL(n, {\Bbb H})$-manifold, $\dim_\R M =4n$,
and 
\[ 
  {\cal R}_{p,q}:\; \Lambda^{p+q,0}_I(M)\arrow \Lambda^{p,q}_{I,+}(M)
\]
the isomorphism induced by $\goth{su}(2)$-action 
as in \ref{_bico_ide_Theorem_}.
Consider the projection 
\begin{equation}\label{_proj_to_+_Equation_} 
\Lambda^{p,q}_{I}(M)\arrow
\Lambda^{p,q}_{I,+}(M),
\end{equation}
and let $R:\; \Lambda^{p,q}_{I}(M)\arrow\Lambda^{p+q,0}_I(M)$
denote the composition of \eqref{_proj_to_+_Equation_} 
and ${\cal R}_{p,q}^{-1}$.

Let $\Phi_I$ be a nowhere degenerate 
holomorphic section of $\Lambda^{2n,0}_I(M)$. Assume that $\Phi_I$ is
real, that is, $J(\Phi_I)=\bar\Phi_I$, and positive.
Existence of such a form is equivalent to 
$\Hol(M) \subset SL(n, {\Bbb H})$ (\ref{_real_holo_parallel_Lemma_}).
It is often convenient to define $SL(n, {\Bbb H})$-structure
by fixing the quaternionic action and the holomorphic
form $\Phi_I$.

\hfill

Define the map
\[ {\cal V}_{p,q}:\;
  \Lambda^{p+q,0}_I(M)\arrow\Lambda^{n+p, n+q}_I(M)
\]
by the relation
\begin{equation}\label{_V_p,q_via_test_form_Equation_}
{\cal V}_{p,q}(\eta) \wedge \alpha = \eta \wedge R(\alpha)\wedge \bar\Phi_I,
\end{equation}
for any test form $\alpha \in \Lambda^{n-p, n-q}_I(M)$.

\hfill

The map ${\cal V}_{p,p}$ is especially remarkable,
because it maps closed, positive
$(2p,0)$-forms to closed, positive $(n+p, n+p)$-forms,
as the following proposition implies.

\hfill

\proposition\label{_V_main_Proposition_}
Let $(M,I,J,K, \Phi_I)$ be an $SL(n, {\Bbb H})$-manifold, and
\[ {\cal V}_{p,q}:\;
  \Lambda^{p+q,0}_I(M)\arrow\Lambda^{4n-p, 4n-q}_I(M)
\]
 the map defined above.
Then
\begin{description}
\item[(i)] ${\cal V}_{p,q}(\eta)= {\cal R}_{p,q}(\eta) \wedge {\cal V}_{0,0}(1)$.
\item[(ii)]  The map ${\cal V}_{p,q}$ is injective, for
  all $p$, $q$.
\item[(iii)] $(\1)^{(n-p)^2}{\cal V}_{p,p}(\eta)$ is real if and
  only $\eta\in\Lambda^{2p,0}_I(M)$ is real, 
and weakly positive if and only if $\eta$ is weakly positive.
\item[(iv)] ${\cal V}_{p,q}(\6\eta)= \6{\cal V}_{p-1,q}(\eta)$,
and ${\cal V}_{p,q}(\6_J\eta)= \bar\6{\cal  V}_{p,q-1}(\eta)$.
\item[(v)] ${\cal V}_{0,0}(1) = \lambda {\cal
  R}_{n,n}(\Phi_I)$, where $\lambda$ is a positive rational number,
depending only on the dimension $n$.
\end{description}

{\bf Proof:} See \cite{_Verbitsky:skoda.tex_}, 
Proposition 4.2, or \cite{_AV:Calabi_}, Theorem 3.6. \endproof

\subsection{Algebra generated by $\omega_I$, $\omega_J$, $\omega_K$}

Let $(M,I,J,K, g)$ be a quaternionic Hermitian manifold.
Consider the algebra $A^*= \oplus A^{2i}$
generated by $\omega_I$, $\omega_J$, and $\omega_K$.
In \cite{_Verbitsky:cohomo_}, this algebra 
was computed explicitly. It was shown that, up to the middle
degree, $A^*$ is a symmetric algebra with generators
$\omega_I$, $\omega_J$, $\omega_K$. The algebra $A^*$ has Hodge
bigrading $A^k = \bigoplus\limits_{p+q=k}A^{p,q}$.
From the Clebsch-Gordan formula, we obtain that
$A^{2i}_+:= \Lambda^{2i}_+(M)\cap A^{2i}$, for $i\leq n$,
is an orthogonal complement to $Q(A^{2i-4})$,
where $Q(\eta) = \eta \wedge (\omega_I^2 +
\omega_J^2+\omega_K^2)$. Moreover, $A^{2i}_+$
is irreducible as a representation of $SU(2)$. 
Therefore, the 
space $A^{p,p}_+= \ker Q^* \restrict {A^{p,p}}$ is 1-dimensional.

\hfill

\proposition\label{_Pi_+_of_omega^n+p_clo_Proposition_}
Let  $(M,I,J,K, \Phi_I)$ be an $SL(n, {\Bbb H})$-manifold,
equipped with an HKT-structure $\Omega$. Assume that
$\Omega^n = \Phi_I$. Let 
\[ \Pi_+:\; \Lambda^{n+q,n+q}_I(M)\arrow \Lambda^{n+q,n+q}_{I,+}(M)\]
be the projection to the component of maximal weight with respect
to the $SU(2)$-action. Then $\Xi_q:=\Pi_+(\omega_I^{n+q,n+q})$ is a
closed, weakly positive $(n+q,n+q)$-form.

\hfill

{\bf Proof:} Consider the algebra $A^*= \oplus A^{2i}$
generated by $\omega_I$, $\omega_J$, and $\omega_K$
as above. The map $R^{p,q}$ is induced by the
$SU(2)$-action, hence it maps $A^{*,*}$ to itself.
Since ${\cal V}_{p,q}(\eta)= 
{\cal R}_{p,q}(\eta) \wedge {\cal V}_{0,0}(1)$,
and ${\cal V}_{0,0}(1)$ is proportional to 
${\cal R}_{n,n}(\Phi_I)\in A^*$, we obtain
\[ {\cal V}_{p,q} (A^{p+q,0}) \subset A^{n+p,n+q}.
\]
Since ${\cal V}(\Omega^p)\subset A^{n+p,n+p}_+$, the 
1-dimensional space $A^{n+p,n+p}_+$ is generated by ${\cal V}(\Omega^p)$.
This form is closed by \ref{_V_main_Proposition_}.
Therefore, the projection of $\omega_I^{n+p}$ to
$A^{n+p,n+p}_+$ is closed.
\endproof


\section{Balanced HKT manifolds}


\subsection{Quaternionic Monge-Amp\`ere equation}
\label{_qCY_Subsection_}

Let $(M,I,J,K)$ be a hypercomplex manifold, and
$\phi \arrow \6\6_J\phi$ the operator
$C^\infty(M) \stackrel{\6\6_J}\arrow \Lambda^{2,0}_I(M)$
defined in Subsection \ref{_posi_2,0-forms_Subsection_}. In \cite{_AV:Calabi_},
the quaternionic Monge-Amp\`ere operator 
$C^\infty(M) \arrow \Lambda^{2n,0}_I(M)$
was defined, mapping a function
$\phi$ to $(\6\6_J\phi)^n$, where $n=\dim_{\Bbb H}M$
(see also \cite{_Alesker:MA_}). This operator is remarkably
similar to the usual Monge-Amp\`ere operators (real and complex)
which are well known in geometry. A quaternionic version
of Calabi-Yau theorem was conjectured.

\hfill

\conjecture\label{_HKT-CY_Conjecture_}
Let $(M,I,J,K, \Omega)$ be a compact HKT manifold
with holonomy $SL(b, {\Bbb H})$,
$\dim_{\Bbb H}M=n$, and $\Phi\in \Lambda^{2n,0}_I(M)$
a nowhere degenerate real\footnote{In the sense of 
Subsection \ref{_posi_2,0-forms_Subsection_}}, 
section of the canonical bundle. Then 
\begin{equation}\label{_q_MA_Equation_}
\Phi=A(\Omega+ \6\6_J\phi)^n
\end{equation}
for some constant $A>0$ and a real function
$\phi\in C^\infty(M)$.

\hfill

It is easy to see that in this case, $\Omega+ \6\6_J\phi$
is an HKT-form (see the proof of \ref{_HKT_CY_Proposition_}),
hence \ref{_HKT-CY_Conjecture_} implies  existence
of an HKT-metric $\Omega'=\Omega+ \6\6_J\phi$
such that the corresponding volume form is proportional
to $\Phi$.

When $\Hol(M) \subset SL(n, {\Bbb H})$, this conjecture
was partly verified in \cite{_AV:Calabi_}.
We have shown that the solution of \eqref{_q_MA_Equation_}
is unique, and also gave a priori $C_0$-bounds on
its solution. For Yau's proof of existence of
solutions of Monge-Amp\`ere equation to work, one
also needs $C_2$ and $C_3$-bounds.

\ref{_HKT-CY_Conjecture_} immediately implies the
following statement. 

\hfill 

\proposition\label{_HKT_CY_Proposition_}
Let $(M,I,J,K, \Omega)$  be a compact HKT-manifold
with holonomy $\Hol(M) \subset SL(n, {\Bbb H})$,
and $\Phi\in \Lambda^{2n,0}_I(M)$ a non-zero section
of the canonical class which is parallel with respect
to the Obata connection $\nabla$. Assume that the
\ref{_HKT-CY_Conjecture_} is true for $M$. Then
$M$ admits an HKT-form $\Omega_1 = \Omega + \6\6_J \phi$
such that $\Omega^n_1$ is a holomorphic volume form.
Moreover, in this case, one has $\nabla(\Omega^n_1)=0$,
where $\nabla$ is the Obata connection.

\hfill

{\bf Proof:} Let $\phi$ be a solution of an equation
$\Phi=A(\Omega+ \6\6_J\phi)^n$.
The form $\Omega_1:= \Omega+ \6\6_J\phi$ is ${\Bbb
  H}$-positive.
Indeed, since $\Phi=A(\Omega+ \6\6_J\phi)^n$,
this form is nowhere degenerate. At a point
$p\in M$ where $\phi$ reaches its minimum, 
the quaternionic Hessian form $\6\6_J\phi$
is positive (Subsection 
\ref{_posi_2,0-forms_Subsection_}), hence at $p$ the 
quaternionic Hermitian form
$x, y \arrow \Omega_1(x, Jy)$ is positive definite.
Since $\Omega_1$ is nowhere degenerate, this form is positive
definite everywhere on $M$. Therefore,
$\Omega_1$ is an HKT-form. 
To check that $\nabla(\Omega_1^n)=0$,
one uses \ref{_real_holo_parallel_Lemma_} below.
We proved
\ref{_HKT_CY_Proposition_}. \endproof

\hfill

The following lemma is essentially contained in 
\cite{_BDV:nilmanifolds_} (Theorem 3.2).

\hfill

\lemma\label{_real_holo_parallel_Lemma_}
Let $(M,I,J,K)$ be a hypercomplex manifold,
and $\eta$ a top degree $(2n, 0)$-form,
which is ${\Bbb H}$-real and holomorphic. Then $\eta$
is Obata-parallel.

\hfill

{\bf Proof:} Since the Obata connection is torsion-free, $d\eta =
\Alt(\nabla \eta)$, where $\Alt= \bigwedge:\; \Lambda^{2n}(M)
\otimes \Lambda^1(M)\arrow\Lambda^{2n+1}(M)$ denotes the exterior
product. Since $\eta$ is holomorphic, $\bar\6\eta=0$. The map $\Alt$
restricted to $\Lambda^{2n,0}(M) \otimes \Lambda^{0,1}(M)$ is an
isomorphism; therefore, $\nabla^{0,1}\eta=0$. Since
$\eta$ is real, $J(\bar\eta)=\eta$, and we have
\[
\nabla^{0,1}J(\bar\eta)=\nabla^{0,1}\eta=0.
\]
This gives $\nabla^{0,1}\bar\eta=0$, because $J$ is Obata-invariant.
 However,
$\nabla^{0,1}\bar\eta=\overline{\nabla^{1,0}\eta}$, and this gives
$\nabla^{1,0}\eta=0$. We proved that $\nabla^{0,1}\eta +
\nabla^{1,0}\eta = \nabla\eta =0$.
\endproof

\hfill

\remark
It is quite hard to construct examples of compact
HKT-manifolds with holonomy in $SL(n, {\Bbb H})$.
So far, the only one construction is known. In
\cite{_BDV:nilmanifolds_}, it was shown that all
hypercomplex nilmanifolds have holonomy in $SL(n, {\Bbb H})$.
However, for an HKT nilmanifold, one always has 
a left-invariant HKT-form (\cite{_Fino_Gra_}),
and such a form satisfies $\nabla(\Omega^n)=0$
(\cite{_BDV:nilmanifolds_}, Theorem 3.2).
Therefore, \ref{_HKT_CY_Proposition_} is true
in this situation, though \ref{_HKT-CY_Conjecture_} 
is not proven even in the simplest cases.

\hfill

\remark \label{_unique_Remark_}
An HKT-form $\Omega_1$ satisfying 
$\Omega_1^n=\phi$ is unique in its cohomology
class $[\Omega_1]\in H^2(\calo_{(M,I)})$,
where $H^2(\calo_{(M,I)})$ denotes the holomorphic
cohomology of $(M,I)$. Indeed, as shown in 
\cite{_Verbitsky:HKT_} (see also \cite{_Alesker_Verbitsky_HKT_}), 
on any $SL(n, {\Bbb H})$-manifold
two closed, real $(2,0)$-forms $\Omega$ and $\Omega'$
which belong to the same cohomology class in $H^2(\calo_{(M,I)})$
satisfy $\Omega-\Omega'=\6\6_J\phi$. However,
the equation \eqref{_q_MA_Equation_} cannot have
two different solutuions, as shown in \cite{_AV:Calabi_},
Theorem 4.7.

\subsection{Balanced HKT-manifolds}

\definition
Let $(M, I, g)$ be a complex Hermitian manifold,
$\dim_\C M =n$, and $\omega\in \Lambda^{1,1}(M)$ its 
Hermitian form. One says that $M$ is {\bf balanced}
if $d (\omega^{n-1})=0$.

\hfill

\remark
It is easy to see that 
$d (\omega^{m})=0$ for $1 \leq m \leq n-2$ implies that
$\omega$ is K\"ahler; the balancedness makes sense
as the only non-trivial condition of form $d (\omega^{m})=0$
which is not equivalent to the K\"ahler property.

\hfill

Let $(M,I,J,K, \Omega)$ be an HKT-manifold,
$\dim_{\Bbb H}M =n$. The ``HKT Calabi-Yau'' condition
$\bar\6(\Omega^n)=0$ has an elegant differential-geometric
interpretation.

\hfill

\theorem \label{_balanced_CY_equiv_Theorem_}
Let $(M,I,J,K, \Omega)$ be an HKT-manifold,
$\dim_{\Bbb H}M =n$. Then the following conditions
are equivalent. 
\begin{description}
\item[(i)] $\bar\6(\Omega^n)=0$
\item[(ii)] $\nabla(\Omega^n)=0$, where $\nabla$ is the Obata connection
\item[(iii)] The manifold $(M,I)$ with the induced quaternionic
Hermitian metric is balanced as a Hermitian manifold: 
\[
d(\omega_I^{2n-1})=0.
\]
\end{description}

\remark
A balanced HKT-manifold has holonomy in $SL(n, {\Bbb
  H})$. 
This statement follows immediately from 
the implication (iii) $\Rightarrow$ (ii)
of \ref{_balanced_CY_equiv_Theorem_}.

\hfill

\remark
The condition $\nabla(\Omega^n)=0$
is independent from the choice of a
basis $I,J,K$, $IJ = -JI=K$ of ${\Bbb H}$.
Indeed, suppose that $g\in SU(n)$, and 
$(I_1, J_1, K_1)= (g(I), g(J), g(K))$
is a new basis in ${\Bbb H}$. 
The corresponding HKT-form
$\Omega_1 = \omega_{J_1}+ \1 \omega_{K_1}$
can be expressed as $\Omega_1 = g(\Omega)$,
hence
\[
\nabla(\Omega_1^n)= \nabla(g(\Omega_1^n)) = g(\nabla(\Omega^n))=0.
\]
Therefore, \ref{_balanced_CY_equiv_Theorem_}
leads to the following corollary.

\hfill

\corollary
Let $(M,I,J,K, \Omega)$ be an HKT-manifold,
such that the corresponding complex Hermitian
manifold $(M,I)$ is balanced. Then $(M, I_1)$
is balanced for any complex structrure $I_1$
induced by the quaternions. Moreover,
$(M,I,J,K, \Omega)$ is an $SL({\Bbb H}, n)$-manifold. 
\endproof

\hfill 

{\bf Proof of \ref{_balanced_CY_equiv_Theorem_}:}
The equivalence (i) $\Leftrightarrow$ (ii)
is immediate from \ref{_real_holo_parallel_Lemma_}.
The implication (ii) $\Rightarrow$ (iii) easily
follows from \ref{_V_main_Proposition_}. Indeed, 
$\Omega^{n-1}$ is $\6$- and $\6_J$-closed, hence
${\cal V}(\Omega^{n-1})$ is a closed $(2n-1, 2n-1)$-form.
This form is proportional to $\omega_I^{2n-1}$,
by \ref{_Pi_+_of_omega^n+p_clo_Proposition_}. 

The implication (iii) $\Rightarrow$ (i)
is proven as follows. 

\hfill

{\bf Step 1:} 
The Hermitian 
manifold $(M, I, \omega_I)$ is balanced if and only if
$d^* \omega_I=0$, which is equivalent to
$\6^*\omega_I = \bar\6^*\omega_I=0$. By \ref{_bico_ide_Theorem_},
this is equivalent to $\6^* \Omega = \6^*_J\Omega=0$.
This gives $d^*\Omega=0$.
We obtain that $d^* \omega_I=0$ leads to 
$d^*\omega_J = d^* \omega_K=0$. This also brings
\[
d\left(\Omega^n \wedge \bar\Omega^{n-1}\right) = *d^*\Omega=0.
\]

{\bf Step 2:}
Let $\bar\theta$ be a $(0,1)$-form defined by
\[
\bar\6\left(\Omega^n\right)= \Omega^n\wedge \bar\theta.
\]
Using a (2,0)-version of Lefschetz
$\goth{sl}(2)$-action (\cite{_Verbitsky:HKT_}), 
it is easy to observe that the map
\begin{equation}\label{_Lef_Omega_Equation_}
\Lambda^{0,1}_I(M) \xlongrightarrow{\wedge \bar\Omega^{n-1}}  \Lambda^{0,n-1}_I(M)
\end{equation}
is an isomorphism. The HKT-condition
gives $\bar\6\bar\Omega=0$, hence
\[
0= \bar\6\left(\Omega^n \wedge \bar\Omega^{n-1}\right)=
\bar\6\left(\Omega^n\right) \wedge \bar\Omega^{n-1}=
\Omega^n\wedge \bar\theta\wedge \bar\Omega^{n-1}.
\]
Using the isomorphism \eqref{_Lef_Omega_Equation_},
we obtain that $\Omega^n\wedge \bar\theta\wedge\bar\Omega^{n-1}=0$
implies that $\bar\theta=0$. Therefore, 
$\bar\6\left(\Omega^n\right)=0$. 
The implication (iii) $\Rightarrow$ (i)
is proven. We have finished the proof of 
\ref{_balanced_CY_equiv_Theorem_}. \endproof

\hfill

\remark
From \ref{_unique_Remark_}, 
we obtain that a balanced HKT-metric on an
HKT-manifold is unique in
its cohomology class 
$[\Omega]\in H^2\left(\calo_{(M,I)}\right)$. In particular,
the space of balanced HKT-metrics is finite-dimensional.

\hfill

\proposition\label{_Xi_k_defi_Proposition_}
Let $(M, I, J, K, \Omega)$ be a balanced HKT manifold,
and $\Xi_k:= \Pi_+(\omega_I^k)$ the
maximal weight component of $\omega_I^k$. Then
$d \Xi_k=0$ for any $k \geq n$, where $n=\dim_{\Bbb H}M$. Moreover,
the $(k,k)$-form $\Xi_k$ is weakly positive.

\hfill

{\bf Proof:} On a balanced HKT-manifold $(M, I, J, K, \Omega)$,
the top exterior power $\Omega^n$ is Obata parallel and
holomorphic (\ref{_balanced_CY_equiv_Theorem_}). Now,
\ref{_Xi_k_defi_Proposition_} is directly implied by
\ref{_Pi_+_of_omega^n+p_clo_Proposition_}. \endproof


\section{Strong HKT manifolds}


\subsection{Strong HKT metrics and HKT-potential}

There is another important class of HKT metrics,
called {\em strong HKT metrics}. For physicists,
such metrics are of special interest (\cite{_GHR_},
\cite{_Howe_Papado_}). The original definition
of HKT metrics (\cite{_Howe_Papado_}) assumed the strong
HKT condition; in mathematical literature it was dropped,
because of relative lack of examples.

\hfill

\definition
Let $(M,I,J,K)$ be a hypercomplex manifold, and
$g$ a quaternionic Hermitian metric on $M$.
The metric $g$ is called {\bf strong HKT} (abbreviated for sHKT)
if it is HKT and, moreover,
\begin{equation}
dd^c\omega_I=0,
\end{equation}
where $d^c:=-IdI$ is the usual twisted differential.

\hfill

\remark\label{_sHKT_Bismut_Remark_}
Let $(M,I,J,K, g)$ be an HKT-manifold.
The torsion of the Bismut connection on $M$ can be written
as
\[
T=-Id\omega_I = -Jd\omega_J = -Kd \omega_K
\]
(\ref{_HKT_via_torsion_Remark_}).
Clearly, the sHKT  condition
is equivalent to $dT=0$.

\hfill

\remark\label{_H_sHKT_Remark_}
Let $(M,I,J,K)$ be a hypercomplex
manifold equipped with a quaternionic Hermitian structure.
Denote by $H$ the $SU(2)$-representation 
generated by the 3-forms $d\omega_I$, $d\omega_J$, $d\omega_K$. 
The HKT condition is equivalent to $\dim H\leq 4$ 
(\ref{_HKT_via_torsion_Remark_}). 
On top of it, the 
strong HKT condition means that for any $v\in H$,
one has $dv=0$. Indeed,
$g$ is HKT if and only if 
\[
Id\omega_I = Jd\omega_J = Kd \omega_K.
\]
Therefore, $H$ is generated by four 3-forms
\begin{equation} \label{_basis_in_H_Equation_}
H=\langle d\omega_I, \ \ I d \omega_I, \ \ J d\omega_I, \ \ K d \omega_I\rangle.
\end{equation}
Writing $d \omega_I = - K d\omega_J= J d \omega_K$,
we obtain $J d \omega_I=-d \omega_K$, $K d \omega_I=d \omega_J$.
This allows us to rewrite the basis \eqref{_basis_in_H_Equation_},
for any HKT-manifold:
\[
H=\langle d\omega_I, \ \ I d \omega_I, 
\ \ J d\omega_I=-d \omega_K, \ \ K d \omega_I=d \omega_J\rangle.
\]
Of those 4 3-forms, 3 are manifestly exact, and 
$I d \omega_I$ is closed if and only if
$M$ is sHKT, as follows from \ref{_sHKT_Bismut_Remark_}.

\hfill

\proposition\label{_sHKT_via_pote_Proposition_}
Let $(M,I,J,K)$ be an HKT-manifold, and $\Omega\in\Lambda^{2,0}_I(M)$
its HKT form. Then the 
strong HKT condition is equivalent to $\bar\6\bar\6_J \Omega=0$. 

\hfill

{\bf Proof:}
Clearly, $\bar\6_J \Omega$ is an element in $H$,
and for all $v\in H$,  and  all differentials $\delta$ of form
$d, IdI, JdJ, KdK$, we have $\delta(v)=0$ per
\ref{_H_sHKT_Remark_}. Therefore, strong HKT implies
$\bar\6\bar\6_J \Omega=0$. 

To prove the converse
implication, consider a local potential $\phi$ for
$\Omega$, $\Omega = \6\6_J \phi$, where $\phi \in C^\infty(M)$
is a smooth function, defined locally on $M$. Such
a potential exists as shown in \cite{_Banos_Swann_} (see also
\cite{_Alesker_Verbitsky_HKT_}). 
Then $\bar\6\bar\6_J \Omega=0$ is equivalent to
\begin{equation}\label{_4-th-order_on_phi_Equation_}
\bar\6\bar\6_J\6\6_J \phi=0.
\end{equation}
The operators $\6$, $\6_J$,
$\bar\6$, $\bar\6_J$ can be written 
down as linear combinations of
twisted differentials $d, d_I:= - I dI, d_J:= -JdJ, d_K:= - KdK$.
These differentials pairwise anticommute
(\cite{_Verbitsky:HKT_}), hence the only
non-zero homogeneous fourth order monomial 
expressed through $d, d_I, d_J, d_K$ is
the operator $\Lambda^*(M) \xlongrightarrow{d d_I d_J  d_K}\Lambda^{*+4}(M)$.
This gives
\begin{equation}\label{_4_th_order_ide_Equation_}
\bar\6\bar\6_J\6\6_J = \const \cdot d d_I d_J  d_K,
\end{equation}
where $\const$ is a non-zero constant which can be found
by an explicit calculation. The Hermitian form $\omega_I$
can be expressed through $\phi$ as
\begin{equation}\label{_HKT_potential_Equation_}
\omega_I = dd_I\phi + d_J d_K \phi
\end{equation}
(\cite{_Banos_Swann_}).
Therefore, the strong HKT condition is equivalent to
\[
0 = dd_I \omega = dd_I d_J d_K \phi=0 
\]
This is equivalent to $\bar\6\bar\6_J\6\6_J\phi=0 $
as follows from \eqref{_4_th_order_ide_Equation_}.
We proved \ref{_sHKT_via_pote_Proposition_}. \endproof

\hfill

\remark 
The strong HKT condition can be expressed through the potential,
as indicated above. An HKT metric with potential
$\phi$ is strong HKT if and only if $\phi$ satisfies
$dd_I d_J d_K \phi=0$.

\subsection{HKT classes and strong HKT metrics}

Let $M$ be a hypercomplex manifold.
Consider the complex
\begin{equation}\label{_HKT_cohomo_complex_}
C^\infty(M) \stackrel{\6\6_J} 
\arrow \Lambda^{2,0}_I(M) \stackrel{\6\oplus \6_J} \arrow \Lambda^{3,0}_I(M) 
\oplus \Lambda^{3,0}_I(M).
\end{equation}
It is easy to see that this complex is elliptic
(\cite{_Alesker_Verbitsky_HKT_}).

Denote the cohomology of the complex 
\eqref{_HKT_cohomo_complex_} by $H^{2,0}_{\6\6_J}(M)$.
From \ref{_bico_ide_Theorem_} it follows
immediately that the group $H^{2,0}_{\6\6_J}(M)$ is 
independent from the choice of a quaternionic triple
$I,J,K$. 

This group is similar to the Bott-Chern cohomology
group in complex geometry
(\cite{_Teleman:cone_}, \cite{_OV:MN_}).
The Bott-Chern cohomology encodes the
 cohomological information about
holomorphic line bundles, and the 
group $H^{2,0}_{\6\6_J}(M)$ encodes the cohomological
information about the HKT-forms.

\hfill

\definition
Let $(M,I,J,K, \Omega)$ be an HKT manifold,
and $[\Omega]\in H^{2,0}_{\6\6_J}(M)$
the cohomology class of the HKT form $\Omega$.
Then $[\Omega]$ is called {\bf the HKT class of $M$}.

Clearly,  HKT forms  $\Omega, \Omega_1$ have the
same HKT class if and only if $\Omega-\Omega_1 = \6\6_J\phi$,
for a smooth function $\phi$, globally defined on $M$.

\hfill

\remark
Using Hodge theory (in particular, the $\6\bar\6$-lemma)
one can prove that the Bott-Chern cohomology
of a compact K\"ahler manifold $X$ is  equal to the usual
Hodge cohomology group $H^{1,1}(X)$. For compact 
HKT-manifolds with holonomy in $SL(n, {\Bbb H})$
a version of Hodge theory was proven in 
\cite{_Verbitsky:HKT_}. In particular,
it was shown that for any compact HKT-manifold with
$\Hol(M) \subset SL(n, {\Bbb  H})$, the natural
map $H^{2,0}_{\6\6_J}(M)\arrow H^2_{\6}(\Lambda^{*,0}(M))\cong \overline{H^2_{\6}(\calo_{(M,I)})}$
to the complex conjugate to the
corresponding holomorphic cohomology is an isomorphism.

\hfill

\remark\label{_overdete_Remark_}
Let $(M,I,J,K)$ be a hypercomplex manifold, 
$\dim_{\Bbb H} M =n$, and
$g$ a quaternionic Hermitian form (such a form
always exists). Denote by $\omega_I$,
$\omega_J$ and $\omega_K$ the corresponding
Hermitian forms \eqref{_three_2_forms_qH_Equation_}. 
Consider the 4-th order operator from $C^\infty(M)$ to $C^\infty(M)$,
\[ 
  \phi 
  \stackrel \Box \arrow 
  (d d_I d_J d_K \phi, \omega_I^2 + \omega_J^2 + \omega_K^2),
\]
where $(\cdot, \cdot)$ denotes the Hermitian product 
on differential forms induced from $g$. Using the usual
K\"ahler identities, one immediately obtains that on
a hyperk\"ahler manifold, $\Box = \Delta^2$, where
$\Delta$ is the Laplacian. Therefore, $\Box$
is elliptic, and the 
equation $d d_I d_J d_K \phi =0$, $\phi \in C^\infty(M)$
is overdetermined.

This leads to the following conjecture.

\hfill

\conjecture\label{_sHKT_unique_Conjecture_}
Let $(M,I,J,K)$ be a compact HKT manifold, and
$[\Omega]\in H^{2,0}_{\6\6_J}(M)$
its HKT class. Then there exists at most
one strong HKT form $\Omega'$ with
$[\Omega']=[\Omega]$.

\hfill

\remark
For $M$ hyperk\"ahler, \ref{_sHKT_unique_Conjecture_}
is clear. Indeed, if strong HKT forms
$\Omega, \Omega'$ have the same 
cohomology class, one has $\Omega-\Omega'=\6\6_J\phi$,
and (as shown in the proof of \ref{_sHKT_via_pote_Proposition_}),
this brings $dd_Id_Jd_K \phi=0$. Using
\ref{_overdete_Remark_}, one obtains that
$\Delta^2\phi=0$, hence $\phi$ is constant.
In particular, any strong HKT metric on 
a compact hyperk\"ahler manifold is automatically
hyperk\"ahler.

\subsection{Strong HKT metrics on balanced HKT manifolds}

A complex Hermitian manifold $(M,I, \omega)$ 
is called {\bf strong KT} if it satisfies $dd^c\omega=0$.
It is balanced if $d^*\omega=0$.
In \cite{_Alexandrov_Ivanov_}, Section 2, Remark 1,
and independently in \cite{_Fino_Salamon_Parton:SKT_}, Proposition 1.4,
it was shown that a compact Hermitian manifold
which is balanced can never be strong KT,
unless it is K\"ahler. A stronger result can
be proven for balanced HKT-manifolds, 
following the same lines, if one uses
the $Sp(1,1)$-representation theory 
as indicated below.

\hfill

\theorem\label{_sHKT_never_bala_Theorem_}
 Let $(M,I,J,K,\Omega)$ be a balanced HKT-manifold.
Assume that $n=\dim_{\Bbb H} M\geq 3$. If $M$ admits a strong HKT-form in
the same cohomology class $[\Omega]\in H^2(\calo_{(M,I)})$, then
$d\Omega=0$, and $M$ is hyperk\"ahler.

\hfill

{\bf Proof:} The balanced HKT condition
is equivalent to $\Lambda_{\omega_I} d\omega_I=0$,
where $\Lambda_{\omega_I}$ is the Hermitian adjoint
to $ \eta \arrow \eta \wedge \omega_I$.
Indeed, 
\[
\Lambda_{\omega_I} d\omega_I = \frac{* d(\omega_I^{2n-1})}{2n-1}.
\]
Forms satisfying $\Lambda_{\omega_I}\eta=0$ are called
{\bf primitive}.

Denote by $\Xi_3\in\Lambda^{2n-3,2n-3}_I(M)$ the maximum weight component of 
$\omega_I^{2n-3}$. By \ref{_Xi_k_defi_Proposition_}, $\Xi_3$ is a closed,
positive $(2n-3, 2n-3)$-form.

The $(2,1)$-Hodge component $(d\omega_I)^{2,1}\in\Lambda^{2,1}_I(M)$
is closed by \ref{_H_sHKT_Remark_},
primitive as seen above, and has weight 1 with respect to
$SU(2)$ as follows from \ref{_bico_ide_Theorem_}. 
Therefore, $(d\omega_I)^{2,1}$
is a highest weight vector with respect to $\goth{so}(1,4)$-action
associated with the quaternionic Hermitian structure
(\cite{_so(5)_}; see Subsection \ref{_so_1,4_Subsection_} for details). 
By \ref{_Hodge_Riemm_hk_Theorem_}, the Hermitian form 
\[
\eta_1, \eta_2 \arrow \1\int_M \eta_1 \wedge \bar\eta_2 \wedge \Xi_3
\]
is sign-definite on the space of primitive $(2,1)$-forms of weight 1.
Therefore,
\[
\int_M (d\omega_I)^{2,1} \wedge (d\omega_I)^{1,2}  \wedge \Xi_3\neq 0
\]
unless $(d\omega_I)^{2,1}=0$.
A trivial calculation gives 
\[ 
(d\omega_I)^{2,1} \wedge (d\omega_I)^{1,2} = \1
d\omega_I\wedge d_I \omega_I. 
\]
Since
\[
0 = \int_M d d_I(\omega_I\wedge \omega_I\wedge \Xi_3)= 
\int_M d\omega_I\wedge d_I \omega_I\wedge \Xi_3 + 
\int_M d d_I\omega_I\wedge \omega_I\wedge \Xi_3
\]
we conclude that
\begin{equation}\label{_dd_I_omega_positi_Equation_}
\int_M d d_I\omega_I\wedge \omega_I\wedge \Xi_3\neq 0,
\end{equation}
unless $(M,I,J,K,\Omega$ is hyperk\"ahler.

Now assume that $(M,I,J,K)$ admits an sHKT-form
$\Omega'$ in the same cohomology class. Then 
$\Omega-\Omega'=\6\6_J\phi$ for some globally defined 
function $\phi$ on $M$. Since
$\bar\6\bar\6_J\Omega'=0$, 
\[
\bar\6\bar\6_J\6\6_J\phi= \bar\6\bar\6_J(\Omega-\Omega')= \bar\6\bar\6_J\Omega.
\]
From  \eqref{_4_th_order_ide_Equation_},
 we obtain that 
\[ 
dd_Id_Jd_K\phi= 
\const\bar\6\bar\6_J\6\6_J\phi= \const\bar\6\bar\6_J\Omega= d d_I\omega_I,
\]
where $\const$ is a non-zero rational number.
Then \eqref{_dd_I_omega_positi_Equation_} gives
\[
0\neq \int_M d d_I\omega_I\wedge \omega_I\wedge \Xi_3=
\int_M  dd_Id_Jd_K\phi\wedge \omega_I\wedge \Xi_3.
\]
Since $\Xi_3$ is $d,d_I,d_J,d_K$-closed, the
 last integral can be rewritten as
\[
\int_M  dd_Id_Jd_K\phi\wedge \omega_I\wedge \Xi_3=
\int_M \phi dd_Id_Jd_K \omega_I\wedge \Xi_3.
\]
However, locally $\omega_I$ is written
as and $\omega_I= dd_I\psi + d_J d_K \psi$
by Banos-Swann \eqref{_HKT_potential_Equation_},
hence the form $dd_Id_Jd_K \omega_I$ is identically zero.
Therefore, the integral
\eqref{_dd_I_omega_positi_Equation_} must vanish as well,
and $(M,I,J,K,\Omega)$ must be hyperk\"ahler.
\endproof

\hfill

\remark
If \ref{_HKT-CY_Conjecture_} 
is true, then every $SL(n, {\Bbb H})$-manifold
contains a balanced HKT-metric in each HKT-class,
necessarily unique. Then \ref{_sHKT_never_bala_Theorem_}
implies that an $SL(n, {\Bbb H})$-manifold admits
no sHKT-metrics, for $n \geq 3$.

\subsection{Strong HKT metrics on nilmanifolds}

The only examples (so far) of compact 
$SL(n, {\Bbb H})$-manifolds are hypercomplex 
nilmanifolds. A hypercomplex nilmanifold
is a quotient of a nilpotent Lie group equipped with
a left-invariant hypercomplex structure, by a
co-compact, discrete subgroup, acting from the left.
In \cite{_Fino_Gra_}  it was shown that any 
HKT-metric on a hypercomplex nilmanifold
can be averaged with the Lie group action,
giving rise to a left-invariant HKT-metric.
The left-invariant HKT-structures can
be considered as metrics on the corresponding
Lie algebra, and studied algebraically.

The HKT-metrics on
hypercomplex nilmanifolds are studied in
\cite{_BDV:nilmanifolds_}. 
It was shown that a
hypercomplex nilmanifold which admits an HKT metric
is necessarily {\bf abelian}, that is, the corresponding
Lie subalgebra ${\goth g}^{1,0}$ of left-invariant 
(1,0)-vector fields is abelian.  In
\cite{_Dotti_Fino:HKT_}, Proposition 2.1
it was shown that any abelian hypercomplex nilmanifold
admitting a strong HKT metric is necessarily
a torus. Therefore, for nilmanifolds 
\ref{_sHKT_never_bala_Theorem_} is known
from conjunction of \cite{_BDV:nilmanifolds_} and
\cite{_Dotti_Fino:HKT_}.

The problem of constructing strong HKT-manifolds
seems to be difficult. All compact non-hyperk\"ahler
examples of sHKT-manifolds known so far are homogeneous.
The hypercomplex structures on semisimple Lie groups
obtained by physicists Ph. Spindel et al (\cite{_SSTvP_}) 
and independently by D. Joyce (\cite{_Joyce_})
are strong HKT (\cite{_Gra_Poon_}). 
A powerful new method of ``doubling'' a
strong HKT $4n$-dimensional Lie algebra to obtain a
strong HKT Lie algebra of dimension $8n$
is proposed in \cite{_Barberis_Fino:sHKT_}.


\section{Appendix: Hyperk\"ahler Hodge-Riemann relations}


\subsection{$\goth{so}(4,1)$-action and the Schur's lemma}
\label{_so_1,4_Subsection_}

Let $V$ be a quaternionic 
Hermitian vector space, and $\Lambda^* V$ its exterior
algebra. Consider the 2-forms $\omega_I$, $\omega_J$, $\omega_K$,
defined on $V$ as in \eqref{_three_2_forms_qH_Equation_}.
Let $\goth a$ be the Lie algebra generated by the operators
\[
L_{\omega_I}(\eta) = \eta \wedge\omega_I, \ \
L_{\omega_J}(\eta) = \eta \wedge\omega_J, \ \ 
L_{\omega_K}(\eta) = \eta \wedge\omega_K,
\]
and their Hermitian adjoints. In \cite{_so(5)_},
it was shown that ${\goth a} \cong \goth{so}(4,1) = \goth{sp}(1,1)$
(in \cite{_V:Mirror_}, the corresponding Lie group was found;
it is isomorphic to $Sp(1,1) \cong \Spin(4,1)$).

\hfill

The Lie algebra $\goth{sp}(1,1)$ has rank 2.
As shown in \cite{_so(5)_}, one could choose a
Cartan subalgebra of $\goth a$ in such a way that
the corresponding weight decomposition of $\Lambda^* V$
coincides with the Hodge decomposition
$\Lambda^* V = \bigoplus_{p,q}\Lambda^{p,q}(V)$.

\hfill

\definition
Let $A$ be a  vector space, and $\goth g$ a Lie
algebra acting on $A$. Assume that $A$ is {\bf a
semisimple representation}, that is, $A$ is a direct
sum of irreducible $\goth g$-representations. 
Consider a decomposition $A= \bigoplus A_\alpha$
of $A$ into a direct sum of representations
of $\goth g$, with each $A_\alpha$ being a sum of isomorphic
irreducible representations, non-isomorphic between
different $A_\alpha$. Such a decomposition is called
{\bf the isotypic decomposition}; it is
obviously unique and well-defined.
\hfill

Let $\Lambda^* V= \oplus_\alpha I_\alpha$ 
be the isotypic decomposition of $\Lambda^* V$
with respect to the action of $\goth a$.
Since the Hodge decomposition on $\Lambda^* V$
is induced by the $\goth{sp}(1,1)$-action, the Hodge
decomposition of  $\Lambda^* V$ is compatible
with the isotypic decomposition. 
Let $I_\alpha= \oplus_{p,q} I^{p,q}_{\alpha}$
be the Hodge decomposition of $I_\alpha$,
taken with respect to the complex structure $I$
on $V$. 

The main result of this Appendix is the following
quaternionic Hermitian version of Hodge-Riemann
relations.

\hfill

\theorem\label{_Hodge_Riemm_hk_Theorem_}
Let $V={\Bbb H}^n$ be a quaternionic Hermitian space,
and $\Lambda^*V= \oplus I^{p,q}_\alpha$ the Hodge decomposition
of its isotypic decomposition defined above. Consider
a form $P\in \Lambda^{2n-k, 2n-k}V$, $P=P(\omega_I, \omega_J, \omega_K)$
obtained as an order $(2n-k)$ homogeneous polynomial 
of $\omega_I, \omega_J, \omega_K$, and let
$(\cdot,\cdot)_P$ be a semi-linear pairing on $I^{p,q}_\alpha$,
defined as
\[
(\eta, \eta')_P:= \frac{\eta\wedge\bar\eta'\wedge P}{\Vol V},
\]
where $p+q=k$ and $\Vol V$ is the Riemannian volume form on $V$. Then
$(\cdot,\cdot)_P$ is sign-definite or identically zero. 

\hfill

\ref{_Hodge_Riemm_hk_Theorem_}
is an immediate consequence of Schur's lemma together
with the following theorem, proven in Subsection \ref{_Howe_dua_Subsection_}.

\hfill

\theorem\label{_Hodge_isoty_irre_Theorem_}
Let $V={\Bbb H}^n$ be a quaternionic Hermitian space,
and $\Lambda^*V= \oplus I^{p,q}_\alpha$ the Hodge decomposition
of its isotypic decomposition defined above. Denote by $Sp(V)$
the group $Sp(n)$ of quaternionic unitary matrices acting on $V$; 
this group obviously commutes with $Sp(1,1)$-action, hence
preserves the spaces $I^{p,q}_\alpha$.
Then the spaces $I^{p,q}_\alpha$ are irreducible as representations
of $Sp(V)$, for all $p,q, \alpha$. 

\hfill

\remark 
Let $W=\C^n$ be a Hermitian vector space,
$\Lambda^*V$ its (real) exterior algebra, equipped
with a usual Lefschetz-type $\goth{sl}(2)$-action,
$\Lambda^*V= \oplus I_\alpha$ its isotypic decomposition,
and $\Lambda^*V= \oplus I^{p,q}_\alpha$ the Hodge
decomposition of $\oplus I_\alpha$.
The Hodge-Riemann relations state that the form
\[
(\eta, \eta'):= \frac{\eta\wedge\bar\eta'\wedge \omega^k}{\Vol V},
\]
is sign-definite on the space $I_\alpha^{p,q}$, where 
$p+q+2k=n$. It is deduced directly from Schur's lemma,
because, as follows from Howe's duality, the spaces $I_\alpha^{p,q}$
are irreducible as representations of $U(n)$
(\cite{_Howe:duality_}), and the form $(\cdot, \cdot)$ is 
$U(n)$-invariant. 

The hyperk\"ahler Hodge-Riemann relations are proven
using the same argument, with $Sp(n)$ instead of $U(n)$.

\hfill

\remark
On a compact hyperk\"ahler manifold $M$, 
 the $Sp(1,1)$-action preserves the harmonic forms.
Therefore, the decomposition $\Lambda^*M = \oplus I^{p,q}_\alpha$
is well defined on harmonic forms. The forms 
$P=P(\omega_I, \omega_J, \omega_K)$ which
can be expressed polynomially through
$\omega_I, \omega_J, \omega_K$ are closed,
hence the pairing $(\cdot,\cdot)$ is
well defined in cohomology. In this situation, the
Hodge-Riemann relations have topological interpretation,
similar to the Hodge index theorem in the K\"ahler case.

\subsection{Howe's duality and $\goth{sp}(1,1)$-action}
\label{_Howe_dua_Subsection_}

Howe's duality can be stated as in 
R. Howe's paper \cite{_Howe:duality_}
in a very general fashion involving
graded Clifford algebras associated with graded
vector spaces. This version of Howe's duality
includes both the usual Clifford algebra and
usual spinors and
its odd counterpart, the Weil algebra (the
algebra of differential operators) and
the Weil representation, also known as 
the space of symplectic spinors.

To obtain the hyperk\"ahler Hodge-Riemann
relations, the symplectic spinorial part of this
picture is not needed.

To simplify the exposition, we omit the odd
Clifford part of the statement, and state
the Howe's duality for usual Clifford algebras
and the usual spinors. 

\hfill

Let $V$ be a vector space, $W= V \oplus V \oplus V \oplus ...$
a sum of several copies of $V$, and $\Lambda^*W$ the 
corresponding Grassmann algebra. Denote by $\tilde W$
the sum $W \oplus W^*$ equipped with a natural symmetric pairing. 
The corresponding Clifford algebra $\Cl(\tilde W)$
is naturaly identified with $\End(\Lambda^*(W))$,
and $\Lambda^*(W)$ is identified with the 
associated spinor space of $\tilde W$.

\hfill

\definition
In these assumptions, let $\Gamma, \Gamma'\subset \goth{o}(\tilde W)$
be Lie subalgebras of the orthogonal Lie algebra $\goth{o}(\tilde W)$.
The pair $(\Gamma, \Gamma')$ is called {\bf a dual pair}
({\bf a Howe's dual pair})
if $\Gamma$ is a centralizer of $\Gamma'$ and
$\Gamma'$ is a centralizer of $\Gamma$.
If $\Gamma$ and $\Gamma'$ are reductive
Lie algebras, $(\Gamma, \Gamma')$ is called 
{\bf a reductive dual pair}.

\hfill

The most important example of a dual pair is provided by the following
general construction.

\hfill

\proposition
(\cite{_Howe:duality_}) Let $G$ be a classical Lie group, $V$ its fundamental
representation, $W= V \oplus V \oplus V \oplus ...$,
and $\tilde W = W \oplus W^*$. Denote by $\Gamma \subset \goth{o}(\tilde W)$
the Lie algebra of $G$ acting on $\tilde W$, and let 
$\Gamma'\subset \goth{o}(\tilde W)$
be its centralizer. Then $(\Gamma, \Gamma')$ is a reductive 
dual pair. \endproof

\hfill

Such a dual pair is called {\bf a classical dual pair}.

The Lie algebra $\goth{o}(\tilde W)$
acts on the corresponding spinor space,
hence we can consider $\goth{o}(\tilde W)$
as a Lie subalgebra in $\Cl(\tilde W)$.

The main result of Howe's duality is the following
useful theorem.

\hfill

\theorem
(Howe's duality)
Let $G$ be a classical Lie group, $V$ its fundamental representation,
and $W= V \oplus V \oplus V \oplus ...$. Consider
the corresponding classical dual pair 
$(\Gamma, \Gamma')\subset \goth{o}(\tilde W)\subset \Cl(\tilde W)$.
Then the associative subalgebra of $\Cl(\tilde W) = \End(\Lambda^*(W))$
generated by $\Gamma'$ is the full algebra of invariants of $G$ in
$\Cl(\tilde W)$. 

{\bf Proof:} This is \cite{_Howe:duality_}, Theorem 7. \endproof

\hfill

To prove \ref{_Hodge_isoty_irre_Theorem_},
let's apply Howe's duality to $G= Sp(V)$,
$W=V$. The centralizer of $G$ in $\goth{o}(\tilde W)$
is naturally identified with $\goth u(1,1, {\Bbb H})= \goth{sp}(1,1)$,
hence $\goth{sp}(1,1)$ generates the full algebra of invariants
of $Sp(V)$ acting on $\Lambda^*(V)$. Since the Cartan
algebra action of $\goth{sp}(1,1)$ coincides
with the Hodge decomposition of $\Lambda^*(V)$, the corresponding
eigenspaces $I^{p,q}_\alpha$ have no endoporphisms
generated by $\goth{sp}(1,1)$. By Howe's duality, this
implies that all $I^{p,q}_\alpha$
are irreducible representations of $Sp(V)$. 
We proved \ref{_Hodge_isoty_irre_Theorem_}.

\hfill

{\bf Acknowledgements:}
I am grateful to Ruxandra Moraru and Semyon Alesker
for their interest and valuable discussions, and to Victor Ginzburg,
Dmitri Kaledin and David Kazhdan for an advice on
representation theory. Special thanks to Isabel Dotti,
Geo Grantcharov and Stefan Ivanov and for insightful 
comments and reference.

\hfill

\hfill

\small{

\noindent {\sc Misha Verbitsky\\
{\sc  Institute of Theoretical and
Experimental Physics \\
B. Cheremushkinskaya, 25, Moscow, 117259, Russia }\\
\tt verbit@mccme.ru }

}


{\small
\begin{thebibliography}{AV1}

\bibitem[AI]{_Alexandrov_Ivanov_}
Bogdan Alexandrov, Stefan Ivanov,
{\em Vanishing theorems on Hermitian manifolds},
arXiv:math/9901090, Diff. Geom. Appl. vol. 14 (2001),
251-263.

\bibitem[A]{_Alesker:MA_} 
Semyon Alesker
{\em Quaternionic Monge-Amp\`ere equations}, arXiv:math/0208005,
	J. Geom. Anal. 13 (2003), no. 2, 205--238.

\bibitem[AV1]{_Alesker_Verbitsky_HKT_} 
Semyon Alesker, Misha Verbitsky,
{\em Plurisubharmonic functions on hypercomplex manifolds and
 HKT-geometry,} arXiv:math/0510140,
J. Geom. Anal. 16 (2006), no. 3, 375--399.



\bibitem[AV2]{_AV:Calabi_} 
Semyon Alesker, Misha Verbitsky
{\em Quaternionic Monge-Amp\`ere equation and Calabi problem 
for HKT-manifolds}, arXiv:0802.4202, 32 pages.

\bibitem[BS]{_Banos_Swann_}
 Banos, Bertrand; Swann, Andrew;
 {\em Potentials for hyper-K\"ahler metrics with torsion},
arXiv:math/0402366, 
 Classical Quantum Gravity 21 (2004), no. 13, 3127--3135.


\bibitem[BDV]{_BDV:nilmanifolds_}
Maria L. Barberis, Isabel G. Dotti, Misha Verbitsky,
{\em Canonical bundles of complex nilmanifolds, 
with applications to hypercomplex geometry},
arXiv:0712.3863, 22 pages.

\bibitem[BF]{_Barberis_Fino:sHKT_}
M. L. Barberis, A. Fino
{\em New strong HKT manifolds arising from quaternionic
representations},  arXiv:0805.2335, 18 pages.


\bibitem[Bes]{_Besse:Einst_Manifo_} 
Besse, A., {\em Einstein Manifolds}, Springer-Verlag, New York (1987)

\bibitem[CS]{_Capria-Salamon_} 
Capria, M. M., Salamon, S. M. 
{\it Yang-Mills fields on quaternionic spaces}, 
Nonlinearity {\bf 1} (1988), no. 4, 517--530. 

\bibitem[DF]{_Dotti_Fino:HKT_}
I. Dotti, A. Fino,  {\em Hyperk\"ahler torsion structures invariant
by nilpotent Lie groups,} Class. Quantum Gravity {\bf 19} (2002),
551--562, math.DG/0112166.




\bibitem[FG]{_Fino_Gra_}  
Fino, A.,  Grantcharov, G.,
{\em On some properties of the manifolds with 
skew-symmetric torsion and holonomy SU(n)
and Sp(n)}, math.DG/0302358, 
Adv. Math. 189 (2004), no. 2, 439--450.

\bibitem[FPS]{_Fino_Salamon_Parton:SKT_}
Anna Fino, Maurizio Parton, Simon Salamon
{\em Families of strong KT structures in six dimensions},
math/0209259, Commentarii Mathematici Helvetici (2004), 2, 317-340.

\bibitem[GHR]{_GHR_}
Gates, S. J., Jr.; Hull, C. M.; Ro\v cek, M.,
{\em Twisted multiplets and new supersymmetric nonlinear $\sigma$-models},
Nuclear Phys. B 248 (1984), no. 1, 157--186. 

\bibitem[G]{_Gauduchon:Herm+Dirac_} 
P. Gauduchon,
{\em Hermitian connections and Dirac operators},
Bollettino U. M. I. B {\bf 11} (1997) 257-288. 

\bibitem[GP]{_Gra_Poon_}
Grantcharov, G., Poon, Y. S.,
{\em Geometry of hyper-K\"ahler connections with torsion},
 math.DG/9908015, 
Comm. Math. Phys. 213 (2000), no. 1, 19--37. 

\bibitem[J]{_Joyce_}
D. Joyce,  {\em Compact hypercomplex and quaternionic manifolds}, J.
Differential Geom. {\bf 35} (1992) no. 3, 743-761



\bibitem[HP]{_Howe_Papado_}
P.S. Howe, G. Papadopoulos,  {\em Twistor spaces for hyper-K\"ahler
manifolds with torsion} Phys. Lett. B 379 (1996), no. 1-4, 80--86.

\bibitem[Ho]{_Howe:duality_}
Howe, Roger, {\em Remarks on classical invariant theory} Trans.
   Amer. Math. Soc. 313 (1989), no. 2, 539--570.

\bibitem[Hu]{_Humphreys_}
Humphreys, J., {\em Introduction to Lie Algebras and Representation
  Theory},  Graduate Texts
in Mathematics, Springer-Verlag, no. 9, 1972.

\bibitem[MS]{_Merkulov_Sch:long_}
S. Merkulov, L. Schwachh\"ofer, {\em  Classification of irreducible
holonomies of torsion-free affine connections}, math.DG/9907206,
Ann. of Math. (2) 150 (1999), no. 1, 77-149, also see Addendum:
math.DG/9911266, Ann. of Math. (2) 150 (1999), no. 3, 1177-1179



\bibitem[Ob]{_Obata_} 
Obata, M., {\em Affine connections on manifolds
with almost complex, quaternionic or Hermitian structure}, 
Jap. J. Math., 26 (1955), 43-79.

\bibitem[OV]{_OV:MN_}
Liviu Ornea, Misha Verbitsky,
{\em Morse-Novikov cohomology of locally conformally K\"ahler manifolds},
arXiv:0712.0107, 22 pages.

\bibitem[SSTV] {_SSTvP_}
Ph. Spindel, A. Sevrin, W. Troost, A. Van Proeyen {\em Extended
supersymmetric $\sigma$-models on group manifolds}, Nucl. Phys. B308
(1988) 662-698.


\bibitem[SSTV]{_Strominger:Bismut_}
A. Strominger, {\em Superstrings with torsion,}
Nuclear Phys. B 274 (1986), no. 2, 253--284.


\bibitem[Te]{_Teleman:cone_}
A. Teleman, {\em
The pseudo-effective cone of a
non-K\"ahlerian surface and applications},
Math. Ann. 335 (2006), no. 4, 965--989.

\bibitem[V0]{_so(5)_} 
M. Verbitsky, 
{\em On the action of the Lie algebra 
$\frak{s}\frak{o}(5)$ on the cohomology
of a hyperk\"ahler manifold}, 
Func. Anal. and Appl. {\bf 24} (1990), 70-71.

\bibitem[V1]{_Verbitsky:cohomo_}
Verbitsky, M., 
{\em  Cohomology of compact hyperk\"ahler manifolds},
  arXiv:alg-geom/9501001,  87 pages.

\bibitem[V2]{_V:Mirror_}
Verbitsky, M.,
{\em Mirror Symmetry for hyperk\"ahler manifolds,}
 alg-geom/9512195,  Mirror symmetry, III (Montreal, PQ, 1995), 115--156,
   AMS/IP Stud. Adv. Math., 10, Amer. Math. Soc., Providence, RI, 1999.



\bibitem[V3]{_Verbitsky:HKT_}
Verbitsky, M., 
{\em 
Hyperk\"ahler manifolds with torsion, supersymmetry and Hodge theory},
math.AG/0112215, Asian J. Math. Vol. 6, No. 4, pp. 679-712 (2002).

\bibitem[V4]{_V:reflexive_}
Verbitsky, M.,
{\em Hyperholomorpic connections   on coherent  sheaves 
and stability}, 40 pages, math.AG/0107182



\bibitem[V5]{_Verbitsky_HKT_exa_}
Verbitsky, M., 
{\em Hyperk\"ahler manifolds with torsion obtained from
hyperholomorphic bundles,}, 
Math. Res. Lett. 10 (2003), no. 4, 501--513, 
also in math.DG/0303129.

\bibitem[V6]{_Verbitsky:canoni_}
M. Verbitsky,  {\em Hypercomplex manifolds with trivial canonical
bundle
 and their holonomy}, arXiv:math/0406537,
``Moscow Seminar on Mathematical Physics, II'',
American Mathematical Society Translations,
{\bf 2}, 221 (2007).

\bibitem[V7]{_Verbitsky:qD_}
M. Verbitsky,
{\em Quaternionic Dolbeault complex and vanishing theorems on
hyperkahler manifolds}, math/0604303, 
Compos. Math. 143 (2007), no. 6, 1576--1592.


\bibitem[V8]{_Verbitsky:skoda.tex_}
Verbitsky, M., {\em Positive forms on hyperkahler
  manifolds}, arXiv:0801.1899, 33 pages.



\end{thebibliography}
}
\end{document}